\documentclass{article}
\usepackage[utf8]{inputenc}
\usepackage{amsfonts}
\usepackage{amssymb}
\usepackage{hyperref}
\usepackage{indentfirst}
\usepackage{mathrsfs}
\usepackage{tikz}
\usepackage{mathtools}
\usepackage{graphicx}
\usepackage{xparse}
\usepackage[all]{xy}

\usepackage[top=1in,bottom=1in,left=1.25in,right=1.25in]{geometry}
\date{}
\allowdisplaybreaks

\begin{document}

 \renewcommand{\baselinestretch}{1.2}
 \renewcommand{\arraystretch}{1.0}

 \title{\bf Braided crossed category over crossed group-cograded weak Hopf quasigroups}
 \author
 { \bf{Huili Liu} \footnote{College of Science, Nanjing Agricultural University, Nanjing 210095, Jiangsu, China.
  E-mail: 2020111006@njau.edu.cn}, \,
  \bf{Lingli Zhu} \footnote{College of Science, Nanjing Agricultural University, Nanjing 210095, Jiangsu, China. E-mail: 2021111005@stu.njau.edu.cn}, \,
  \bf{Tao Yang} \footnote{Corresponding author. College of Science, Nanjing Agricultural University, Nanjing 210095, Jiangsu, China.
             E-mail: tao.yang@njau.edu.cn}
 }
 \maketitle

 \begin{center}
 \begin{minipage}{12.cm}
 {Abstract:In this paper, we generalizing the main result in Liu\cite{LYZ22} to weak Hopf coquasigroups case. We first define and study group-cograded weak Hopf quasigroups,  which generalize both  group-cograded Hopf quasigroups and weak Hopf group-coalgebras. Then we introduce the notion of $p$-Yetter-Drinfeld weak quasimodule over group-cograded weak Hopf quasigroups $H$. If the antipode of $H$ is bijective,  we show that the category $\mathcal Y\mathcal D\mathcal W\mathcal Q(H)$ of Yetter-Drinfeld weak quasimodules over $H$ is a crossed category, and the subcategory $\mathcal Y\mathcal D(H)$ of Yetter-Drinfeld modules is a braided crossed category.
 }
 \\

 { Key words: weak Hopf quasigroup; crossed group-cograded weak Hopf quasigroup; p-Yetter-Drinfeld weak quasimodule; braided crossed category}
 \\

 { Mathematics Subject Classification 2020: 16T05, 17A01, 18M15}
 \end{minipage}
 \end{center}
 \normalsize

 \section{Introduction}
 \def\theequation{\thesection.\arabic{equation}}
 \setcounter{equation}{0}

 Yetter-Drinfeld modules over bialgebras were introduced by Yetter in \cite{Y90}, also known as Yang-Baxter modules in \cite{LR93} and Yetter-Drinfeld structures in \cite{RT93}.
 The Yetter-Drinfeld modules over bialgebras can form a pre-braided monoidal category which is the main feature of this definition.
 If endow bialgebras with an antipode to obtain Hopf algebras, then the category of Yetter-Drinfeld modules over Hopf algebras is braided.
 The (pre)-braiding structures of Yetter-Drinfeld modules plays a part in the relations between knot theory and quantum.
 This concept was soon generalized to other algebraic structures such as weak Hopf algebras\cite{CWY05}, Hom-bialgebras \cite{MP14}, Hopf group-coalgebras\cite{Z04}, multiplier Hopf group-coalgebras\cite{D08} et al.

 The notion of weak Hopf quasigroups was introduced by Alvarez et al. \cite{AFG16}. As a new generalization of Hopf algebras, weak Hopf quasigroups encompass weak Hopf algebras and Hopf quasigroups, under the unified approach, the more relevant properties of these algebraic structures can be captured.
 The first non-trivial examples of this algebraic structures can be obtained by considering bigroupoids, that is considering the bicategories where every 1-cell is an equivalence and every 2-cell is an isomorphism.
 Since weak Hopf quasigroups satisfy a big group of common properties with weak Hopf algebras, $\rm\acute{A}$lvarez et al. in \cite{AFG21} recently obtain a categorical equivalence between the category of pointed cosemisimple weak Hopf quasigroups over $K$ and a generalization of the category of finite groupoids to the non associative setting.

 It is have been proved that the categories of Yetter-Drinfeld modules over weak Hopf quasigroups\cite{YW21} and group-cograded Hopf algebras are braided \cite{D08}.
 Then there has a natural question: if we endow the weak Hopf quasigroups with a group-cograded structure, is it has the similar results with the weak Hopf quasigroups and group-cograded Hopf algebras case?
 Following this idea, we expanded our study of this article.

 This paper is organized as follows: In Section 2, we recall some notions, such as braided crossed categories, Turaev’s left index notation, and weak Hopf quasigroups.
 All of these are the most important building blocks for completing this article.

 In Section 3, we give the notion of group-cograded weak Hopf quasigroup and introduce some properties of group-cograded weak Hopf quasigroup. Moreover, we define a crossed structure on group-cograded weak Hopf quasigroup, and give a method to obtain a new crossed group-cograded weak Hopf quasigroup.

 In Section 4, we introduce the definition of $p$-Yetter-Drinfeld weak quasimodules over $H$, then show the category $\mathcal Y\mathcal D\mathcal W\mathcal Q(H)$ of Yetter-Drinfeld weak quasimodules is a crossed category, and the subcategory $\mathcal Y\mathcal D(H)$ of Yetter-Drinfeld modules over $H$ is a braided crossed category.

%

 \section{Preliminaries}

 \subsection{Crossed categories}

 Recall from \cite{JS93, ML71, T} that
 let $G$ be a group and $Aut(\mathscr C)$ be the group of invertible strict tensor functors from $\mathscr C$ to itself. A category $\mathscr C$ over $G$ is called a crossed category if it satisfies the following:
 \begin{itemize}
   \item[(1)] $\mathscr C$ is a monoidal category;
   \item[(2)] $\mathscr C$ is disjoint union of a family of subcategories $(\mathscr C_{\alpha})_{\alpha\in G}$, and for any $U\in \mathscr C_{\alpha}, V\in \mathscr C_{\beta}, U\otimes V\in \mathscr C_{\alpha \beta}$. The subcategory $\mathscr C_{\alpha}$ is called the $\alpha th$ component of $\mathscr C$;
   \item[(3)] Consider a group homomorphism $\phi : G\rightarrow Aut(\mathscr C), \beta \rightarrow \phi_{\beta}$, and assume that $\phi_{\beta}(\mathscr C_{\alpha})=\mathscr C_{\beta \alpha\beta^{-1}}$, for all $\alpha,\beta\in G$. The functors $\phi_{\beta}$ are called conjugation isomorphisms.
 \end{itemize}

 We will use Turave's left index notation from \cite{T} for functors $\phi_{\beta}$: Given $\beta \in G$ and an object $V\in \mathscr C$, the functor $\phi_{\beta}$ will be denoted by $^{\beta}(\cdot)$ or $^{V}(\cdot)$ and $^{\beta^{-1}}(\cdot)$ will be denoted by ${}^{\bar{V}}(\cdot)$. Since ${}^{V}(\cdot)$ is a functor, for any object $U\in \mathscr C$ and any composition of morphism $g\circ f$ in $\mathscr C$, we obtain ${}^{V}id_{U}=id_{{}^{V}U}$ and $^{V}(g\circ f)={}^{V}g\circ {}^{V}f$. Since the conjugation $\varphi: G\rightarrow Aut(\mathscr C)$ is a group homomorphism, for any $V, W\in \mathscr C$, we have ${}^{V\otimes W}(\cdot)={}^{V}({}^{W}(\cdot)$ and $^{1}(\cdot)={}^{V}(^{\bar{V}}(\cdot))={}^{\bar{V}}({}^{V}(\cdot))=id_{\mathscr C}$. Since for any $V\in \mathscr C$, the functor $^{V}(\cdot)$ is strict, we have ${}^{V}(g\otimes f)={}^{V}g\otimes {}^{V}f$ for any morphism $f$ and $g$ in $\mathscr C$, and ${}^{V}(1)=1$.
 \\

 A braiding of a crossed category $\mathscr C$ is a family of isomorphisms $(C= C_{U,V})_{U,V\in \mathscr C}$, where
 $C_{U,V}: U\otimes V\rightarrow {}^{U}V\otimes U$satisfying the following conditions:

 \begin{itemize}
   \item[(1)] For any arrow $f\in \mathscr C_{p}(U,U')$ and $g\in \mathscr C(V,V')$,
   \begin{eqnarray}
   (({}^{p}g)\otimes f)\, C_{U,V}=C_{U',V'}\, (f\otimes g);
   \end{eqnarray}
   \item[(2)] For all $U, V, W\in \mathscr C$, we have
   \begin{eqnarray}
   C_{U\otimes V , W}&=&a_{{}^{U\otimes V}W, U, V}\, (C_{U,{}^{V}W}\otimes id_{V})\, a^{-1}_{U, {}^{V}W, V}\, (\iota_{U}\otimes C_{V,W})\, a_{U,V,W},\\
   C_{U, V\otimes W}&=&a^{-1}_{{}^{U}V, {}^{U}W, U}\, (\iota_{{}^{U}V}\otimes C_{U, W})\, a_{{}^{U}V, U, W}\, (C_{U,V}\otimes \iota_{W})\, a^{-1}_{U,V,W},
   \end{eqnarray}
   where a is the natural isomorphisms in the tensor category $\mathscr C$.
   \item[(3)] For all $U, V\in \mathscr C$ and $q\in G$,
   \begin{eqnarray}
   \phi_{q}(C_{U,V})=C_{\phi_{q}(U),\phi_{q}(V)}.
   \end{eqnarray}
 \end{itemize}

 A crossed category endowed with a braiding is called a braided crossed category.
 \\

 Recall from \cite{K95} we can assigned any strict tensor category $\mathcal C$ to a braided tensor category $\mathcal Z(\mathcal C)$. An object of $\mathcal Z(\mathcal C)$ is a pair $(V, C_{-,V})$ where $V$ is an object of $\mathcal C$ and $C_{-,V}$ is a family of natural isomorphisms
 \begin{eqnarray}
 && C_{X,V}: X\otimes V\rightarrow sV\otimes X
 \end{eqnarray}
 defined for all objects $X\in \mathcal C$ such that for all objects $X,Y \in \mathcal C$ we have
 \begin{eqnarray}
 && C_{X\otimes Y,V}=(C_{X,V}\otimes id_{Y})(id_{X}\otimes C_{Y,X}).
 \end{eqnarray}

 A morphism from $(V,C_{-,V})$ to $(W,C_{-,W})$ is a morphism $f: V\rightarrow W$ in $\mathcal C$ such that for each object $X$ of $\mathcal C$ we have
 \begin{eqnarray}
 && (f\otimes id_{X})C_{X,V}=C_{X,W}(id_{X}\otimes f).
 \end{eqnarray}

 The naturality of $C_{-,V}$ means that the following diagram
  $$\xymatrix{
  X\otimes V \ar[d]_{f\otimes id_{V}} \ar[r]^{C_{X,V}}
                & V\otimes X \ar[d]^{id_{V}\otimes f}  \\
  Y\otimes V  \ar[r]_{C_{Y,V}}
                & V\otimes Y             }$$
 commutes for any morphism.

 According to the definition of $\mathcal Z(\mathcal C)$ we have the identity $id_{V}$ is a morphisms in $\mathcal Z(\mathcal C)$ and that if $f,g$ are composable morphisms in $\mathcal Z(\mathcal C)$ then the composition $g\circ f$ in $\mathcal C$ is a morphism in $\mathcal Z(\mathcal C)$. Consequently, $\mathcal Z(\mathcal C)$ is a category in  which the identity of $(V, C_{-,V})$ is $id_{V}$.


\subsection{Weak Hopf quasigroups}

 Throughout this article, all spaces we considered are over a fixed field $k$.

 Recall from \cite{AFG16} that a weak Hopf quasigroup $H$ is a unital magma $(H, \mu, \eta$ and a comonoid $(H, \delta, \epsilon)$ such that the following axioms hold: for any $h,g,l\in H$,
 \begin{itemize}
   \item [(1)]$\delta(hg)=h_{1}g_{1}\otimes h_{2}g_{2}.$

   \item [(2)]$\epsilon((hg)l)=\epsilon(h(gl))=\epsilon(hg_{1})\epsilon(g_{2}l)=\epsilon(hg_{2})\epsilon(g_{1}l).$

   \item [(3)]$1_{1}\otimes 1_{2}\otimes 1_{3}=1_{1}\otimes 1_{2}1'_{1}\otimes 1'_{2}=1_{1}\otimes 1'_{1}1_{2}\otimes 1'_{2}.$

   \item [(4)] There exists a linear map $S: H\rightarrow H$ (called the antipode of $H$) such that, if we define the target morphism $\epsilon_{t}: H\rightarrow H$ by $\epsilon_{t}=id_{H}*S$ and the source morphism $\epsilon_{s}: H\rightarrow H$ by $\epsilon_{s}=S*id_{H}$, for any $h,g \in H$,
       \begin{itemize}
         \item [(4-1)] $\epsilon_{t}(h)=\epsilon(1_{1}h)1_{2};$
         \item [(4-2)]$\epsilon_{s}(h)=1_{1}\epsilon(h1_{2});$
         \item [(4-3)]$S=S*\epsilon_{t}=\epsilon_{s}*S;$
         \item [(4-4)]$S(h_{1})(h_{2}g)=\epsilon_{s}(h)g;$
         \item [(4-5)]$h_{1}(S(h_{2}g))=\epsilon_{t}(h)g;$
         \item [(4-6)]$(hg_{1})S(g_{2})=h\epsilon_{t}(g);$
         \item [(4-7)]$(hS(g_{1}))g_{2}=h\epsilon_{s}(g).$
       \end{itemize}
 \end{itemize}

 If $H$ is a weak Hopf quasigroup, the equalities
 \begin{eqnarray}
 &&\epsilon_{t}*id_{H}=id_{H}*\epsilon_{s}=id_{H},\\
 &&\epsilon_{t}\eta=\eta=\epsilon_{s}\eta,\\
 &&\epsilon\epsilon_{t}=\epsilon=\epsilon\epsilon_{s}
 \end{eqnarray}
 are satisfied, the antipode of weak Hopf quasigroup $H$ is unique and $S(1)=1, \epsilon S=\epsilon$ and for any $h,g\in H$,
 \begin{eqnarray}
 &&S(hg)=S(g)S(h),\\
 &&S(h)_{1}\otimes S(h)_{2}=S(h_{2})\otimes S(h_{1}).
 \end{eqnarray}

 More details of weak Hopf quasigroup see \cite{AFG15}.

%

\section{Crossed group-cograded weak Hopf quasigroup}
 \def\theequation{\thesection.\arabic{equation}}
 \setcounter{equation}{0}

 Motivated by the Hopf quasigroups case, we give the notion of group-cograded weak Hopf quasigroup in this section, and introduce some properties of  group-cograded weak Hopf quasigroup. Moreover, we give a crossed structure on group-cograded weak Hopf quasigroups, and we can use the mirror reflection to obtain a new crossed group-cograded weak Hopf quasigroup.
 \\

 {\bf Definition \thesection.1} A  group-cograded weak Hopf quasigroup $H=(H_{p}, \mu_{p}, \eta_{p}, \Delta,\epsilon)_{p\in G}$ is a family of algebras $(H_{p}, \mu_{p}, \eta_{p})$, may not be associative, and a group coalgebra $(H_{p}, \Delta=(\Delta_{p, q}), \epsilon=\epsilon_{e})_{p,q\in G}$ such that the following conditions hold:
 \begin{itemize}
   \item [(1)] The comultiplication $\Delta_{p,q}: H_{pq}\rightarrow H_{p}\otimes H_{q}$ is a homomorphism of algebras such that
    \begin{eqnarray}
    &&(\Delta_{p,q}\otimes id_{H_{r}})\Delta_{pq,r}(1_{pqr})=(\Delta_{p,q}(1_{pq})\otimes 1_{r})(1_{p}\otimes \Delta_{q,r}(1_{qr})),\\
    &&(\Delta_{p,q}\otimes id_{H_{r}})\Delta_{pq,r}(1_{pqr})=(1_{p}\otimes \Delta_{q,r}(1_{qr}))(\Delta_{p,q}(1_{pq})\otimes 1_{r}),
    \end{eqnarray}
    for all $p, q, r\in G$.

   \item [(2)] The counit $\epsilon: H_{e} \rightarrow k$ is a $k$-linear map satisfying the identity
    \begin{eqnarray}
    \epsilon(ghl)=\epsilon(gh_{(2,e)}\epsilon(h_{(1, e)}l)=\epsilon(gh_{(1,e)})\epsilon(h_{(2,e)}l),
    \end{eqnarray}
    for all $g, h , l\in H_{e}$.

   \item [(3)] There exists a family of $k$-linear maps $S=(S_{p}: H_{p}\rightarrow H_{p^{-1}})_{p\in G}$ such that, if a family of linear maps $\epsilon_{t}=(\epsilon^{t}_{p}: H_{e}\rightarrow H_{p})_{p\in G}$ and $\epsilon_{s}=(\epsilon^{s}_{p}: H_{e}\rightarrow H_{p})_{p\in G}$ defined by
       \begin{eqnarray}
        \epsilon^{t}_{p}(h)=\epsilon(1_{(1,e)}h)1_{(2,p)},\\
        \epsilon^{s}_{p}(h)=1_{(1,p)}\epsilon(h1_{(2, e)}),
       \end{eqnarray}
       for any $h\in H_{e}$, where $\epsilon_{t}$ and $\epsilon_{s}$ are called the $G$-target and $G$-source counit morphisms, for any $h\in H_{e}, g\in H_{p}$,
       \begin{eqnarray}
        &&S=S*\epsilon_{t}=\epsilon_{s}*S,\\
        &&S(h_{(1,p^{-1})})\big(h_{(2,p)}g\big)=\epsilon^{s}_{p}(h)g,\label{SE1}\\
        &&h_{(1,p)}\big(S(h_{(2,p^{-1})}g)\big)=\epsilon^{t}_{p}(h)g,\label{SE2}\\
        &&\big(gh_{(1,p)}\big)S(h_{(2,p^{-1})})=g\epsilon^{t}_{p}(h),\label{SE3}\\
        &&\big(gS(h_{(1,p^{-1})})\big)h_{(2,p)}=g\epsilon^{s}_{p}(h).\label{SE4}
       \end{eqnarray}
 \end{itemize}

 Let $H$ be a group-cograded weak Hopf quasigroup. We define $\tilde{\epsilon}_{t}=(\tilde{\epsilon}^{t}_{p}: H_{e}\rightarrow H_{p})_{p\in G}$ and $\tilde{\epsilon}_{s}=(\tilde{\epsilon}^{s}_{p}: H_{e}\rightarrow H_{p})_{p\in G}$ by the following: for all $p\in G, h\in H_{e}$,
 \begin{eqnarray}
 \tilde{\epsilon}^{t}_{p}(h)=1_{(1,p)}\epsilon(1_{(2,e)}h),\\
 \tilde{\epsilon}^{s}_{p}(h)=\epsilon(h 1_{(1, e)})1_{(2,p)}.
 \end{eqnarray}



 It is easy to check that for all $p\in G$, $\epsilon_{p}^{t}, \epsilon_{p}^{s}, \tilde{\epsilon}_{p}^{t}$ and $\tilde{\epsilon}_{p}^{s}$ are idempotent, and it is possible to prove the following conditions which involve these four morphisms and the antipode S(see \cite{AFG16}): for all $p\in G$,
 \begin{eqnarray}
 &&\epsilon^{t}_{p}\tilde{\epsilon}^{t}_{p}=\epsilon^{t}_{p},\quad \epsilon^{t}_{p}\tilde{\epsilon}^{s}_{p}=\tilde{\epsilon}^{s}_{p},\\
 &&\tilde{\epsilon}^{t}_{p}\epsilon^{t}_{p}=\tilde{\epsilon}^{t}_{p},\quad \tilde{\epsilon}^{s}_{p}\epsilon^{t}_{p}=\epsilon^{t}_{p},\\
 &&\epsilon^{s}_{p}\tilde{\epsilon}^{t}_{p}=\tilde{\epsilon}^{t}_{p},\quad \epsilon^{s}_{p}\tilde{\epsilon}^{s}_{p}=\epsilon^{s}_{p},\\
 &&\tilde{\epsilon}^{t}_{p}\epsilon^{s}_{p}=\epsilon^{s}_{p},\quad \tilde{\epsilon}^{s}_{p}\epsilon^{s}_{p}=\tilde{\epsilon}^{s}_{p},\\
 &&\epsilon^{t}_{p}=\tilde{\epsilon}^{s}_{p}S_{p^{-1}}=S_{p^{-1}}\tilde{\epsilon}^{t}_{p^{-1}},\quad \epsilon^{s}_{p}=\tilde{\epsilon}^{t}_{p}S_{p^{-1}}=S_{p^{-1}}\tilde{\epsilon}^{s}_{p^{-1}},\\
 &&\epsilon^{t}_{p}S_{p^{-1}}=\epsilon^{t}_{p}\epsilon^{s}_{p}=S_{p^{-1}}\epsilon^{s}_{p^{-1}},\quad \epsilon^{s}_{p}S_{p^{-1}}=\epsilon^{s}_{p}\epsilon^{t}_{p}=S_{p^{-1}}\epsilon^{t}_{p^{-1}}.
 \end{eqnarray}

 Some properties of group-cograded weak Hopf quasigroups are given as follows, and all of these properties can be proved with the similar ways in \cite{AFG16,YW21}:
 \begin{itemize}
   \item [(1)] for all $p, q\in G$ and $h\in H_{p}, g\in H_{e}$, we have
   \begin{eqnarray}
   &&h\epsilon^{t}_{p}(g)=\epsilon(h_{(1, e)}g)h_{(2, p)},\quad \epsilon^{s}_{p}(g)h=h_{(1,p)}\epsilon(gh_{(2,e)}),\\
   &&h\tilde{\epsilon}^{t}_{p}(g)=h_{(1,p)}\epsilon(h_{(2,e)}g),\quad \tilde{\epsilon}^{s}_{p}(g)h=\epsilon(gh_{(1,e)})h_{(2,p)},\\
   &&h_{(1,p)}\otimes \epsilon^{t}_{q}(h_{(2,e)})=1_{(1,p)}h\otimes 1_{(2,q)},\label{het1}\\
   &&\epsilon^{s}_{q}(h_{(1,e)})\otimes h_{(2,p)}=1_{(1,q)}\otimes h1_{(2,p)},\label{esh1}\\
   &&\tilde{\epsilon}^{t}_{q}(h_{(1,e)})\otimes h_{(2,p)}=1_{(1,q)}\otimes1_{(2,p)}h,\label{eth1}\\
   &&h_{(1,p)}\otimes \tilde{\epsilon}^{s}_{q}(h_{(2,e)})=h1_{(1,p)}\otimes 1_{(2,q)}\label{hes1},\\
   &&1_{(1,p)}\otimes \epsilon^{t}_{q}(1_{(2,e)})=1_{(1,p)}\otimes 1_{(2,q)},\label{1et}\\
   &&\epsilon^{s}_{q}(1_{(1,e)})\otimes 1_{(2,p)}=1_{(1,q)}\otimes 1_{(2,p)},\label{es1}\\
   &&\tilde{\epsilon}^{t}_{q}(1_{(1,e)})\otimes 1_{(2,p)}=1_{(1,q)}\otimes1_{(2,p)},\label{et1}\\
   &&1_{(1,p)}\otimes \tilde{\epsilon}^{s}_{q}(1_{(2,e)})=1_{(1,p)}\otimes 1_{(2,q)},\label{1es}\\
   &&\epsilon^{t}_{q}(h_{(1,e)})\otimes h_{(2,p)}=S(1_{(1,q^{-1})})\otimes 1_{(2,p)}h,\label{SE5}\\
   &&h_{(1,p)}\otimes \epsilon^{s}_{q}(h_{(2,e)})=h1_{(1,p)}\otimes S(1_{(2,q^{-1})}).\label{SE6}
   \end{eqnarray}

   \item [(2)] for all $h, g\in H_{e}$,
   \begin{eqnarray}
   &&\epsilon\big(h\epsilon^{t}_{e}(g)\big)=\epsilon(hg)=\epsilon\big(\epsilon^{s}_{e}(h)g\big),\label{eets1}\\
   &&\epsilon\big(h\tilde{\epsilon}^{t}_{e}(g)\big)=\epsilon(hg)=\epsilon\big(\tilde{\epsilon}^{s}_{e}(h)g\big),\label{eets2}\\
   &&\epsilon^{t}_{p}\big(h\epsilon_{e}(g)\big)=\epsilon^{t}_{p}(hg),\quad \epsilon^{s}_{p}\big(\epsilon_{e}(h)g\big)=\epsilon^{s}_{p}(hg),\\
   &&\epsilon^{t}_{p}(h)_{(1,p)}\otimes \epsilon^{t}_{q}\big(\epsilon^{t}_{e}(h)_{(2,e)}\big)=\epsilon^{t}_{p}(h)_{(1,p)}\otimes \epsilon^{t}_{q}(h)_{(2,q)},\\
   &&\epsilon^{s}_{p}\big(\epsilon^{s}_{e}(h)_{(1,e)}\big)\otimes \epsilon^{t}_{q}(h)_{(2,q)}=\epsilon^{s}_{p}(h)_{(1,p)}\otimes \epsilon^{t}_{q}(h)_{(2,q)},\\
   &&\epsilon^{t}_{p}\big(\epsilon^{t}_{e}(h)g\big)=\epsilon^{t}_{p}(h)\epsilon^{t}_{e}(g),\quad \epsilon^{s}_{p}\big(h\epsilon^{s}_{e}(g)\big)=\epsilon^{s}_{p}(h)\epsilon^{s}_{e}(g),\\
   &&\tilde{\epsilon}^{t}_{p}\big(\tilde{\epsilon}^{t}_{e}(h)g\big)=\tilde{\epsilon}^{t}_{p}(h)\tilde{\epsilon}^{t}_{e}(g),\quad \tilde{\epsilon}^{s}_{p}\big(h\tilde{\epsilon}^{s}_{e}(g)\big)=\tilde{\epsilon}^{s}_{p}(h)\tilde{\epsilon}^{s}_{e}(g).
   \end{eqnarray}

   \item [(3)] for all $h, g\in H_{p}$,
   \begin{eqnarray}
   gh&=&\big(g\epsilon^{t}_{p}(h_{(1,e)})\big)h_{(2,p)}=\epsilon^{t}_{p}(g_{(1,e)})(g_{(2,p)}h)\\ \nonumber
   &=&g_{(1,p)}\big(\epsilon^{s}_{p}(g_{(2,e)})h\big)=(gh_{(1,p)})\epsilon^{s}_{p}(h_{(2,e)}).
   \end{eqnarray}
 \end{itemize}

 {\bf Definition \thesection.2} A group-cograded weak Hopf quasigroup $H=(H_{p}, \mu_{p}, \eta_{p}, \Delta, \epsilon, S)_{p\in G}$ is called a crossed group-cograded weak Hopf quasigroup if endowed with a family of algebra isomorphisms $\pi=(\pi_{p}: H_{q}\rightarrow H_{pqp^{-1}})_{p,q\in G}$ (called a crossing) such that $(\pi_{p}\otimes \pi_{p})\Delta_{q,r}=\Delta_{pqp^{-1}, prp^{-1}}\pi_{p}$, $\epsilon\pi_{p}=\epsilon$, $\pi_{pq}=\pi_{p}\pi_{q}$, for all $p,q,r\in G$.
 \\

 If $H$ is crossed with the crossing $\pi=(\pi_{p})_{p\in G}$, then we have
 \begin{eqnarray}
 &&\pi_{q}\epsilon^{s}_{p}=\epsilon^{s}_{qpq^{-1}}\pi_{q},\\
 &&\pi_{q}\epsilon^{t}_{p}=\epsilon^{t}_{qpq^{-1}}\pi_{q},
 \end{eqnarray}
 for all $p,q \in G$.
 \\

 {\bf Remark} A (crossed) group-cograded weak Hopf quasigroup is a (crossed) group-cograded Hopf quasigroup if and only if the counit is a homomorphisms of algebras, on the other hand, we know that weak Hopf quasigroup is a generalization of weak Hopf algebras, based on this relationship, we conclude that if the (crossed) group-cograded weak Hopf quasigroup is associative, then it is indeed a (crossed) weak Hopf group-coalgebra.
 \\

 Turaev in \cite{T} proposed a mirror structure for crossed Hopf group-coalgebra, as for a crossed group-cograded weak Hopf quasigroup, we can define the mirror structure in a similar way.

 {\bf Proposition \thesection.3} Let $(H=\bigoplus_{p\in G}H_{p}, \Delta, \epsilon, S, \pi)$ be a crossed group-cograded weak Hopf quasigroup, then we can define its mirror $(\widetilde{H}=\bigoplus_{p\in G}\widetilde{H}_{p}, \widetilde{\Delta}, \widetilde{\epsilon}, \widetilde{S}, \widetilde{\pi})$ by the following way:
 \begin{itemize}
   \item[(1)] as an algebra, $\widetilde{H}_{p}=H_{p^{-1}}$, for all $p\in G$;

   \item[(2)] define the comultiplication $\widetilde{\Delta}_{p,q}: \widetilde{H}_{pq}\rightarrow \widetilde{H}_{p}\otimes \widetilde{H}_{q}$ by : for $h_{q^{-1}p^{-1}}\in \widetilde{H}_{pq}$,
   \begin{eqnarray}
   \widetilde{\Delta}_{p, q}(h_{q^{-1}p^{-1}})&=&(\pi_{q}\otimes id_{H_{q^{-1}}})\Delta_{q^{-1}p^{-1}q, q^{-1}}(h_{q^{-1}p^{-1}});\label{th3.52}
   \end{eqnarray}

   \item[(3)] the counit $\widetilde{\epsilon}$ of $\widetilde{H}$ is the original counit $\epsilon$;

   \item[(4)] the antipode $\widetilde{S}_{p}=\pi_{p}S_{p^{-1}}: \widetilde{H}_{p}=H_{p^{-1}}\rightarrow H_{p}=\widetilde{H}_{p^{-1}}$;

   \item[(5)] for all $p\in G$, define the cross action $\widetilde{\pi}_{p}=\pi_{p}$.
 \end{itemize}
 Then $(\widetilde{H}=\bigoplus_{p\in G}\widetilde{H}_{p}, \widetilde{\Delta}, \widetilde{\epsilon}, \widetilde{S}, \widetilde{\pi})$ is also a crossed group-cograded weak Hopf quasigroup.
 \\

{\bf\emph{Proof}} Straightforward. $\hfill \square$
 \\

 We end this section by some examples of crossed group-cograded weak Hopf quasigroups, and both examples are derived from an action of $G$ on a weak Hopf quasigroup $H, \delta ,\epsilon, S$ over $k$ by weak Hopf quasigroups endomorphisms.

 {\bf Example \thesection.4}  Set $H^{G}=\big(H_{p})_{p\in G}$ and $G$ is the homomorphism group of weak Hopf quasigroup $H$, where for each $p \in G$, the algebra
 $H_{p}$ is a copy of $H$. Fix an identification isomorphism of algebras $i_{p}: H\rightarrow H_{p}$. For $p, q \in G$, we define a comultiplication
 $\Delta_{p, q}: H_{pq}\rightarrow H_{p}\otimes H_{q}$ by
 \begin{eqnarray*}
 &&\Delta_{p, q}(i_{pq}(h))=\sum_{(h)}i_{p}(h_{(1)})\otimes i_{q}(h_{(2)}),
 \end{eqnarray*}
 where $h\in H$. The counit $\epsilon: H_{e}\rightarrow k$ is defined by $\epsilon(i_{e}(h))=\epsilon(h)\in k$ for $h\in H$. For $p \in G$, the antipode $S_{p}: H_{p}\rightarrow H_{p^{-1}}$ is given by
 \begin{eqnarray*}
 && S_{p}(i_{p}(h))=i_{p^{-1}}(S(h)),
 \end{eqnarray*}
 where $h\in H$. For $p, q \in G$, the homomorphism $\pi_{p}: H_{q}\rightarrow H_{pqp^{-1}}$ is defined by $\pi_{p}(i_{q}(h))=i_{pqp^{-1}}(p(h))$. It is easy to check that $H^{G}$ is a crossed group-cograded weak Hopf quasigroup.
 \\

 {\bf Example \thesection.5} We will consider the set $H^{G}$ which is introduced in Example \thesection.3.
 Set $\widetilde{H}^{G}$ be the same family of algebras $(H_{p}=H)_{p \in G}$ with the same counit and
 the same action $\pi$ of $G$, the comultiplication $\widetilde{\Delta}_{p, q}: H_{pq}\rightarrow H_{p}\otimes H_{q}$,
 and the antipode $\widetilde{S}_{p}: H_{p}\rightarrow H_{p^{-1}}$ defined by
 \begin{eqnarray*}
 &&\widetilde{\Delta}_{p, q}(i_{pq}(h))=\sum_{(h)}i_{p}(q(h_{(e)}))\otimes i_{q}(h_{(2)}),\\
 &&\widetilde{S}_{p}(i_{p}(h))=i_{p^{-1}}(p(S(h)))=i_{p^{-1}}(S(p(h))),
 \end{eqnarray*}
 where $h\in H$. The axioms of a crossed group-cograded weak Hopf quasigroup for $\widetilde{H}^{G}$ follow from definitions.
 \\

 {\bf Example \thesection.6} We consider the set $\widetilde{H}^{G}$ as above only with the new comultiplication.
 Set $\bar{H}^{G}=\widetilde{H}^{G}$ with the same counit, the same homomorphism $\pi$ and the same antipode $S$ as ones in $\widetilde{H}^{G}$, but with a new comultiplication
 \begin{eqnarray*}
 &&\bar{\Delta}_{p, q}: H_{pq}\rightarrow H_{p}\otimes H_{q},\qquad \bar{\Delta}_{p, q}(i_{pq}(h))=\sum_{(h)}i_{p}((h_{(1)}))\otimes i_{q}(p(h_{(2)})),\\
 \end{eqnarray*}
 where $h\in H$. The axioms of a (crossed) group-cograded weak Hopf quasigroup for $\bar{H}^{G}$ follow from definitions.
 \\

 Note that $H^{G}$, $\widetilde{H}^{G}$ and $\bar{H}^{G}$ are extensions of $H$ since $H^{G}_{e}=\widetilde{H}^{G}_{e}=\bar{H}^{G}_{e}=H_{e}$
 as weak Hopf quasigroups. The crossed group-cograded weak Hopf quasigroups $H^{G}$ and $\widetilde{H}^{G}$, which defined in Example \thesection.4 and \thesection.5, respectively, are mirrors of each other.
 \\

 \section{ Construction of braided crossed categories}
 \def\theequation{\thesection.\arabic{equation}}
 \setcounter{equation}{0}

 Let $H=\bigoplus_{r\in G}H_{r}$ is a crossed group-cograded weak Hopf quasigroup with a bijective antipode $S$. We introduce the definition of $p$-Yetter-Drinfeld weak quasimodules over $H$, then show the category $\mathcal Y\mathcal D\mathcal W\mathcal Q(H)$ of Yetter-Drinfeld weak quasimodules is a crossed category, and the subcategory $\mathcal Y\mathcal D(H)$ of Yetter-Drinfeld modules over $H$ is a braided crossed category.
 \\

 {\bf Definition \thesection.1}
 Let $H$ be a crossed group-cograded weak Hopf quasigroup over group $G$ and $p$ be a fixed element in $G$. A pair $(V, \rho^{V}=(\rho^{V}_{r})_{r\in G})$ is a right-right $p$-Yetter-Drinfeld weak quasimodule over H, if it satisfies the following conditions:
 \begin{itemize}
   \item[(1)] for all $p\in G$, $V$ is a right weak $H_{p}$-quasimodule in the sense that: for all $h\in H_{e}, v\in V$,
   \begin{eqnarray}
   && v\cdot 1_{p}=v,\nonumber \\
   &&(v\cdot h_{(1,p)})\cdot S_{p^{-1}}(h_{(2,p^{-1})})=v\cdot \epsilon^{t}_{p}(h),\label{WQM1}\\
   &&(v\cdot S_{p^{-1}}(h_{(1,p^{-1})}))\cdot h_{(2,p)}=v\cdot \epsilon^{s}_{p}(h).\label{WQM2}
   \end{eqnarray}

   \item[(2)] for any $r\in G$, $\rho^{V}_{r}: V\rightarrow V\otimes H_{r}$ is a $k$-linear morphism, such that $V$ is coassociativity, that is, for any $r_{1},r_{2}\in G$,
   \begin{eqnarray*}
   (\rho_{r_{1}}^{V}\otimes id_{H_{r_{2}}})\rho_{r_{2}}^{V}=(id_{V}\otimes \Delta_{r_{1},r_{2}})\rho_{r_{1},r_{2}}^{V};
   \end{eqnarray*}
   and $V$ is counitary in the sense that
   \begin{eqnarray*}
   (id_{V}\otimes \epsilon)\rho_{e}^{V}=id_{V};
   \end{eqnarray*}

   \item[(3)] V is crossed in the sense that for all $v\in V$, $r\in G$ and $h,g\in H_{r}$,
   \begin{eqnarray}
   &\rho^{V}_{r}(v\cdot h)=v_{(0)}\cdot h_{(2,p)}\otimes S\pi_{p^{-1}}(h_{(1,pr^{-1}p^{-1})})\big(v_{(1,r)}h_{(3,r)}\big), \label{JR1}\\
   &v_{(0)}\otimes (gk)v_{(1,r)}=v_{(0)}\otimes g(k v_{(1,r)}), \label{JH1}\\
   &v_{(0)}\otimes (gv_{(1,r)})k=v_{(0)}\otimes g(v_{(1,r)}k), \label{JH2}
   \end{eqnarray}
   where $\rho^{V}_{r}(v\cdot h)=(v\cdot h)_{(0)}\otimes (v\cdot h)_{(1,r)}$, for all $h\in H_{p}, g,k\in H_{r}$.
 \end{itemize}

 If $S$ is bijective, for a right weak $H_{p}$-quasimodule $V$, we have the following conditions: for all $h\in H_{e}, v\in V$
 \begin{eqnarray}
 &&(v\cdot S^{-1}(h_{(2,p)})\cdot h_{(1,p^{-1})}=v\cdot \tilde{\epsilon}_{p}^{t}(h),\label{WQM3}\\
 &&(v\cdot h_{(2,p^{-1})})\cdot S^{-1}(h_{(1,p)}=v\cdot \tilde{\epsilon}_{p}^{s}(h).\label{WQM4}
 \end{eqnarray}

 Given two $p$-Yetter-Drinfeld weak quasimodules $(V, \rho{}^{V})$ and $(W, \rho^{W})$, a morphism of this two $p$-Yetter-Drinfeld weak quasimodules $f: (V, \rho^{V})\rightarrow (W,\rho^{W})$ is an $H_{p}$-linear map $f: V\rightarrow W$ such that for any $r\in G$,
 \begin{eqnarray*}
 \rho_{r}^{W}f=(f\otimes id_{H_{r}})\rho_{r}^{V},
 \end{eqnarray*}
 i.e., for all $v\in V$,
 \begin{eqnarray*}
 f(v)_{(0)}\otimes f(v)_{(1,r)}=f(v_{(0)})\otimes v_{(1,r)}.
 \end{eqnarray*}

 Then we have the category $\mathcal Y\mathcal D\mathcal W\mathcal Q(H)_{p}$ of $p$-Yetter-Drinfeld weak quasimodules, the composition of morphisms of $p$-Yetter-Drinfeld weak quasimodules is the standard composition of the underlying linear maps. Moreover, if we assume that $V$ is a right $H_{p}$-module, then we say that is a right-right $p$-Yetter-Drinfeld module. Obviously, right-right $p$-Yetter-Drinfeld modules with the obvious morphisms is a subcategory of $\mathcal Y\mathcal D\mathcal W\mathcal Q(H)_{p}$, and denote it by $\mathcal Y\mathcal D(H)_{p}$.
 \\

 {\bf Proposition \thesection.2} Let $H$ be a crossed group cograded weak Hopf quasigroup. We assume that equality (\ref{JR1}) holds, and $V=(V, \cdot, \rho^{V})$ is a triple where $(V, \cdot)$ is a right weak$H_{p}$-quasimodule, and $(V, \rho^{V})$ is a right $H$-comodule.  The following identities hold: for any $h\in H_{p}, g\in H_{e}$ and $v\in V$, we have
 \begin{eqnarray}
 &&(v\cdot \epsilon^{t}_{e}(g))\cdot h=v\cdot \epsilon^{t}_{e}(g)h,\\
 &&(v\cdot h)\cdot \epsilon^{s}_{e}(g)=v\cdot h\epsilon^{t}_{e}(g),\\
 &&(v\cdot \tilde{\epsilon}^{t}_{e}(g))\cdot h=v\cdot \tilde{\epsilon}^{t}_{e}(g)h,\\
 &&(v\cdot h)\cdot \tilde{\epsilon}^{s}_{e}(g)=v\cdot h\tilde{\epsilon}^{s}_{e}(g).
 \end{eqnarray}
 and
 \begin{eqnarray}
 &&(v\cdot 1_{(2,p)})\cdot h\otimes 1_{(1,q)}= v\cdot 1_{(2,p)}h\otimes 1_{(1,q)},\label{v1h1}\\
 &&(v\cdot h)\cdot 1_{(1,p)}\otimes 1_{(2,q)}= v\cdot h 1_{(1,p)}\otimes 1_{(2,q)},\label{v1h2}\\
 &&(v\cdot 1_{(1,p)})\cdot h\otimes 1_{(2,q)}= v\cdot 1_{(1,p)}h\otimes 1_{(2,q)},\label{v1h3}\\
 &&(v\cdot h)\cdot 1_{(2,p)}\otimes 1_{(1,q)}= v\cdot h 1_{(2,p)}\otimes 1_{(1,q)}.\label{v1h4}
 \end{eqnarray}

{\bf\emph{Proof}} It is straightforward by the similar way in \cite{YW21}. $\hfill \square$.
 \\

 {\bf Proposition \thesection.3} The equation (\ref{JR1}) is equivalent to
 \begin{eqnarray}
 &&v_{(0)}\cdot h_{(1,p)}\otimes v_{(1,r)}h_{(2,r)}=(v\cdot h_{(2,p)})_{(0)}\otimes \pi_{p^{-1}}(h_{(1,prp^{-1})})(v\cdot h_{(2,p)})_{(1,r)}, \label{JR2}\\
 &&v_{(0)}\cdot 1_{(1,p)}\otimes v_{(1,r)}1_{(2,r)}=v_{(0)}\otimes v_{(1,r)}, \label{JR3}
 \end{eqnarray}
 for all $h\in H_{pr}, v\in V$.
 \smallskip

{\bf\emph{Proof}}
 Indeed, suppose the condition (\ref{JR1}) holds, then we have
 \begin{eqnarray*}
 &&(v\cdot h_{(2,p)})_{(0)}\otimes \pi_{p^{-1}}(h_{(1,prp^{-1})})(v\cdot h_{(2,p)})_{(1,r)}\\
 &=&v_{(0)}\cdot h_{(3,p)}\otimes \pi_{p^{-1}}(h_{(1,prp^{-1})})\big(S\pi_{p^{-1}}(h_{(2,pr^{-1}p^{-1})})(v_{(1,r)}h_{(4,r)})\big)\\
 &=&v_{(0)}\cdot h_{(3,p)}\otimes \pi_{p^{-1}}\Big(h_{(1,prp^{-1})}\big(S(h_{(2,pr^{-1}p^{-1})})\pi_{p}(v_{(1,r)}h_{(4,r)})\big)\Big)\\
 &\stackrel{(\ref{SE2})}=&v_{(0)}\cdot h_{(2,p)}\otimes \pi_{p^{-1}}\big(\epsilon^{t}_{prp^{-1}}(h_{(1,e)})\pi_{p}(v_{(1,r)}h_{(3,r)})\big)\\
 &=&v_{(0)}\cdot h_{(2,p)}\otimes \pi_{p^{-1}}\epsilon^{t}_{prp^{-1}}(h_{(1,e)}))(v_{(1,r)}h_{(3,r)})\\
 &\stackrel{(\ref{SE5})}=&v_{(0)}\cdot 1_{(2,p)}h_{(1,p)}\otimes \pi_{p^{-1}}S(1_{(1,pr^{-1}p^{-1})})(v_{(1,r)}h_{(2,r)})\\
 &\stackrel{(\ref{v1h1})}=&(v_{(0)}\cdot 1_{(2,p)})\cdot h_{(1,p)}\otimes \pi_{p^{-1}}S(1_{(1,pr^{-1}p^{-1})})(v_{(1,r)}h_{(2,r)})\\
 &=&(v_{(0)}\cdot 1_{(2,p)})\cdot h_{(1,p)}\otimes \pi_{p^{-1}}S(1_{(1,pr^{-1}p^{-1})})(v_{(1,r)}h_{(2,r)})\\
 &=&(v_{(0)}\cdot1_{(2,p)})\cdot 1'_{(1,p)}h_{(1,p)}\otimes \pi_{p^{-1}}S(1_{(1,pr^{-1}p^{-1})})\big(v_{(1,r)}(1'_{(2,r)}h_{(2,r)})\big)\\
 &=&(v_{(0)}\cdot1_{(2,p)}1'_{(1,p)})\cdot h_{(1,p)}\otimes \pi_{p^{-1}}S(1_{(1,pr^{-1}p^{-1})})\big(v_{(1,r)}(1'_{(2,r)}h_{(2,r)})\big)\\
 &=&(v_{(0)}\cdot1_{(2,p)})\cdot h_{(1,p)}\otimes \pi_{p^{-1}}S(1_{(1,pr^{-1}p^{-1})})\big(v_{(1,r)}(1_{(3,r)}h_{(2,r)})\big)\\
 &\stackrel{(\ref{JH1}),(\ref{JH2})}=&(v_{(0)}\cdot1_{(2,p)})\cdot h_{(1,p)}\otimes \pi_{p^{-1}}\big(S(1_{(1,pr^{-1}p^{-1})})(v_{(1,r)}1_{(3,r)})\big)h_{(2,r)}\\
 &\stackrel{(\ref{JR1})}=&(v\cdot 1)_{(0)}\cdot h_{(1,p)}\otimes (v\cdot 1)_{(1,r)}h_{(2,r)}\\
 &=&v\cdot h_{(1,p)}\otimes v_{(1,r)}h_{(2,r)},
 \end{eqnarray*}
 and
 \begin{eqnarray*}
 &&v_{(0)}\cdot 1_{(1,p)}\otimes v_{(1,r)}1_{(2,r)}\\
 &\stackrel{(\ref{JR2})}=&(v\cdot 1_{(2,p)})_{(0)}\otimes \pi_{p^{-1}}(1_{(1,prp^{-1})})(v\cdot 1_{(2,p)})_{(1,r)}\\
 &\stackrel{(\ref{JR1})}=&v_{(0)}\cdot 1_{(3,p)}\otimes \pi_{p^{-1}}(1_{(1,prp^{-1})})\big(S\pi_{p^{-1}}(1_{(2,pr^{-1}p^{-1})})(v_{(1,r)}1_{(4,r)})\big)\\
 &=&v_{(0)}\cdot 1_{(3,p)}\otimes \pi_{p^{-1}}\Big(1_{(1,prp^{-1})}\big(S(1_{(2,pr^{-1}p^{-1})})\pi_{p}(v_{(1,r)}1_{(4,r)})\big)\Big)\\
 &\stackrel{(\ref{SE2})}=&v_{(0)}\cdot 1_{(2,p)}\otimes \pi_{p^{-1}}\big(\epsilon^{t}_{prp^{-1}}(1_{(1,e)})\pi_{p}(v_{(1,r)}1_{(3,r)})\big)\\
 &=&v_{(0)}\cdot 1_{(2,p)}\otimes \pi_{p^{-1}}\epsilon^{t}_{prp^{-1}}(1_{(1,e)})(v_{(1,r)}1_{(3,r)})\\
 &\stackrel{(\ref{SE5})}=&v_{(0)}\cdot 1_{(2,p)}\otimes \pi_{p^{-1}}S(1_{(1,prp^{-1})})(v_{(1,r)}1_{(3,r)})\\
 &\stackrel{(\ref{JR1})}=&(v\cdot 1)_{(0)}\otimes (v\cdot 1)_{(1,r)}\\
 &=&v_{(0)}\otimes v_{(1,r)},
 \end{eqnarray*}
 thus we conclude that if (\ref{JR1}) holds, (\ref{JR2}) and (\ref{JR3}) hold.

 Conversely, if (\ref{JR2}) and (\ref{JR3}) hold, then
 \begin{eqnarray*}
 &&v_{(0)}\cdot h_{(2,p)}\otimes S\pi_{p^{-1}}(h_{(1,pr^{-1}p^{-1})})(v_{(1,r)}h_{(3,r)}), \\
 &\stackrel{(\ref{JR2})}=&(v\cdot h_{(3,p)})_{(0)}\otimes S\pi_{p^{-1}}(h_{(1,pr^{-1}p^{-1})})\big(\pi_{p^{-1}}(h_{(2,prp^{-1})})(v\cdot h_{(3,p)})_{(1,r)}\big)\\
 &\stackrel{(\ref{SE1})}=&(v\cdot h_{(2,p)})_{(0)}\otimes \pi_{p^{-1}}\epsilon^{s}_{prp^{-1}}(h_{(1,e)})(v\cdot h_{(2,p)})_{(1,r)})\\
 &\stackrel{(\ref{esh1})}=&(v\cdot h_{(2,p)}1_{(2,p)})_{(0)}\otimes \pi_{p^{-1}}(1_{prp^{-1}})(v\cdot h_{(2,p)}1_{(2,p)})_{(1,r)})\\
 &\stackrel{(\ref{v1h4})}=&(v\cdot h_{(2,p)}\cdot 1_{(2,p)})_{(0)}\otimes \pi_{p^{-1}}(1_{prp^{-1}})(v\cdot h_{(2,p)}\cdot 1_{(2,p)})_{(1,r)})\\
 &\stackrel{(\ref{JR2})}=&(v\cdot h)_{(0)}\cdot 1_{(1,p)}\otimes (v\cdot h)_{(1,r)}\cdot 1_{(2,r)}\\
 &\stackrel{(\ref{JR3})}=&(v\cdot h)_{(0)}\otimes (v\cdot h)_{(1,r)},
 \end{eqnarray*}
 hence the equality (\ref{JR1}) holds. $\hfill \square$
 \\

 {\bf Proposition \thesection.4} Let $H$ be a crossed group cograded weak Hopf quasigroup and $V$ is a $p$-Yetter-Drinfeld weak quasimodule (for all $p\in G$) such that the following equations hold:
 \begin{eqnarray}
 &&v_{(0)}\otimes w\cdot h\cdot v_{(1,p)}=v_{(0)}\otimes w\cdot hv_{(1,p)},\\
 &&v_{(0)}\otimes w\cdot v_{(1,p)}\cdot h=v_{(0)}\otimes w\cdot v_{(1,p)}h,
 \end{eqnarray}
 for any $h\in H_{p}, v,w\in V$ and $p\in G$. Then the map
 \begin{eqnarray}
 c: V\otimes V\rightarrow V\otimes V : c(v\otimes w)=w_{(0)}\otimes v\cdot w_{(1,p)},
 \end{eqnarray}
 for any $v,w\in V$ is a solution of the $\mathbf{QYBE}$.

 {\bf\emph{Proof}} It is easy to check that the result is true.$\hfill \square$
 \\

 {\bf Proposition \thesection.5} If $(V, \rho^{V})\in \mathcal Y\mathcal D\mathcal W\mathcal Q(H)_{p}$ and $(W, \rho^{W})\in \mathcal Y\mathcal D\mathcal W\mathcal Q(H)_{q}$, then $V\otimes_{s_{pq}} W\in \mathcal Y\mathcal D\mathcal W\mathcal Q(H)_{pq}$ with the module and comodule structures as follows:
 \begin{eqnarray}
 (v\otimes_{s_{pq}} w)\cdot h_{pq}&=& v\cdot h_{(1,p)}\otimes_{s_{pq}} w\cdot h_{(2,q)}, \label{HHQ1}\\
 \rho_{r}^{V\otimes_{s_{pq}} W}(v\otimes_{s_{pq}} w)&=&v_{(0)}\otimes_{s_{pq}} w_{(0)}\otimes \pi_{q^{-1}}(v_{(1,qrq^{-1})})w_{(1,r)} \label{HHQ2},
 \end{eqnarray}
 where $v\in V, w\in W$ and $h_{pq}\in H_{pq}$.

 {\bf\emph{Proof}} We first check that $V\otimes_{s_{pq}} W$ is a right $H_{pq}$-weak quasimodule, and the unital property is obvious.
 We only check the equation (\ref{WQM1}), the equation (\ref{WQM2})is similar. Indeed , for all $v\otimes w\in V\otimes_{s_{pq}} W$,
 \begin{eqnarray*}
 &&\big((v\otimes_{s_{pq}} w)\cdot h_{(1,pq)}\big)\cdot S(h_{(2,q^{-1}p^{-1})})\\
 &\stackrel{(\ref{HHQ1})}=&(v\cdot h_{(1,p)}\otimes_{s_{pq}} w\cdot h_{(2,q)})\cdot S(h_{(2,q^{-1}p^{-1})})\\
 &\stackrel{(\ref{HHQ2})}=&(v\cdot h_{(1,p)})\cdot S(h_{(4,p^{-1})})\otimes_{s_{pq}} (w\cdot h_{(2,q)})\cdot S(h_{(3,q^{-1})})\\
 &\stackrel{(\ref{WQM1})}=&(v\cdot h_{(1,p)})\cdot S(h_{(3,p^{-1})})\otimes_{s_{pq}} w\cdot \epsilon^{t}_{q}(h_{(2,e)})\\
 &\stackrel{(\ref{het1})}=&(v\cdot 1_{(1,p)}h_{(1,p)})\cdot S(h_{(2,p^{-1})})\otimes_{s_{pq}} w\cdot 1_{(2,q)}\\
 &\stackrel{(\ref{v1h3})}=&(v\cdot 1_{(1,p)})\cdot h_{(1,p)}\cdot S(h_{(2,p^{-1})})\otimes_{s_{pq}} w\cdot 1_{(2,q)}\\
 &\stackrel{(\ref{WQM1})}=&(v\cdot 1_{(1,p)})\cdot \epsilon^{t}_{p}(h)\otimes_{s_{pq}} w\cdot 1_{(2,q)}\\
 &=&(v\cdot 1_{(1,p)})\cdot 1'_{(2,p)}\epsilon(1'_{(1,e)}h)\otimes w\cdot 1_{(2,q)}\\
 &=&v\cdot 1_{(1,p)}1'_{(2,p)}\epsilon(1'_{(1,e)}h)\otimes w\cdot 1_{(2,q)}\\
 &=&v\cdot 1_{(2,p)}\epsilon(1_{(1,e)}h)\otimes w\cdot 1_{(3,q)}\\
 &=&(v\otimes w)\cdot 1_{(2,pq)}\epsilon(1_{(1,e)}h)\\
 &=&(v\otimes_{s_{pq}} w)\cdot \epsilon^{t}_{pq}(h),
 \end{eqnarray*}
 since then $V\otimes_{s_{pq}} W$ is a right $H_{pq}$-weak quasimodule.

 The coassociativity is easy to check, the next we will prove the counitary of $v\otimes_{s_{pq}} w$ hold: for all $v\otimes_{s_{pq}} w\in  V\otimes_{s_{pq}} W$, we have
 \begin{eqnarray*}
 &&(id_{V\otimes_{s_{pq}} W}\otimes \epsilon)\rho^{V\otimes_{s_{pq}} W}(v\otimes_{s_{pq}} w)\\
 &=&v_{(0)}\otimes_{s_{pq}} w_{(0)}\otimes \epsilon(\pi_{q^{-1}}(v_{(1,e)})w_{(1,e)})\\
 &=&v_{(0)}\otimes_{s_{pq}} w_{(0)}\epsilon(\pi_{q^{-1}}(v_{(1,e)})\pi_{q^{-1}}(1')_{1,e})\epsilon(\pi_{q^{-1}}(1')_{(2,e)}w_{(1,e)})\\
 &\stackrel{(\ref{JR3})}=&v_{(0)}\otimes_{s_{pq}} (w_{(0)}\cdot 1_{(1,q)})\epsilon(\pi_{q^{-1}}(v_{(1,e)})\pi_{q^{-1}}(1')_{1,e})\epsilon\big(\pi_{q^{-1}}(1')_{(2,e)}(w_{(1,e)}1_{(2,e)})\big)\\
 &\stackrel{(\ref{JR2})}=&v_{(0)}\otimes_{s_{pq}} (w\cdot 1_{(2,q)})_{(0)}\epsilon\pi_{q^{-1}}(v_{(1,e)}1'_{1,e})\epsilon\big(\pi_{q^{-1}}(1'_{(2,e)})(\pi_{q^{-1}}(1_{(1,e)})(w\cdot 1_{(2,q)})_{(1,e)}\big)\\
 &\stackrel{(\ref{JH1})}=&v_{(0)}\otimes_{s_{pq}} (w\cdot 1_{(2,q)})_{(0)}\epsilon(v_{(1,e)}1'_{1,e})\epsilon\big(\pi_{q^{-1}}(1'_{(2,e)}1_{(1,e)})(w\cdot 1_{(2,q)})_{(1,e)}\\
 &=&v_{(0)}\otimes_{s_{pq}} (w\cdot 1_{(3,q)})_{(0)}\epsilon(v_{(1,e)}1_{1,e})\epsilon\big(\pi_{q^{-1}}(1_{(2,e)})(w\cdot 1_{(3,q)})_{(1,e)}\big)\\
 &\stackrel{(\ref{JR1})}=&v_{(0)}\otimes_{s_{pq}} (w\cdot 1_{(4,q)})\epsilon(v_{(1,e)}1_{(1,e)})\epsilon\Big(\pi_{q^{-1}}(1_{(2,e)})\big(S\pi_{q^{-1}}(1_{(3,e)})(w_{(1,e)}1_{(5,e)})\big)\Big)\\
 &\stackrel{(\ref{SE2})}=&v_{(0)}\otimes_{s_{pq}} (w\cdot 1_{(3,q)})\epsilon(v_{(1,e)}1_{(1,e)})\epsilon\big(\pi_{q^{-1}}\epsilon^{t}_{e(1_{(2,e)})}(w_{(1,e)}1_{(4,e)})\big)\\
 &=&v_{(0)}\otimes_{s_{pq}} (w\cdot 1_{(3,q)})\epsilon(v_{(1,e)}1_{(1,e)})\epsilon\big(\tilde{\epsilon}^{s}_{e}S\pi_{q^{-1}}(1_{(2,e)})(w_{(1,e)}1_{(4,e)})\big)\\
 &\stackrel{(\ref{eets1})}=&v_{(0)}\otimes_{s_{pq}} (w\cdot 1_{(2,q)})_{(0)}\epsilon(v_{(1,e)}1_{(1,e)})\epsilon((w\cdot 1_{(2,q)})_{(1,e)})\\
 &=&(v_{(0)}\cdot 1_{(2,p)})_{(0)}\otimes_{s_{pq}} w\cdot 1_{(3,q)}\epsilon\big(\pi_{p^{-1}}(1_{(1,e)})(v\cdot 1_{(2,p)})_{(1,e)}\big)\\
 &=&v_{(0)}\cdot 1_{(3,p)}\otimes_{s_{pq}}w\cdot 1_{(5,q)}\epsilon\Big(\pi_{p^{-1}}(1_{(1,e)})\big(S\pi_{p^{-1}}(1_{(2,e)})(v_{(1,e)}1_{(4,e)})\big)\Big)\\
 &=&v_{(0)}\cdot 1_{(2,p)}\otimes_{s_{pq}}w\cdot 1_{(4,q)}\epsilon\big(\pi_{p^{-1}}\epsilon^{t}_{e})(1_{(1,e)})(v_{(1,e)}1_{(3,e)})\big)\\
 &=&v_{(0)}\cdot 1_{(2,p)}\otimes_{s_{pq}}w\cdot 1_{(4,q)}\epsilon\big(\tilde{\epsilon}^{s}_{e}S\pi_{p^{-1}}(1_{(1,e)})(v_{(1,e)}1_{(3,e)})\big)\\
 &=&(v\cdot 1_{(1,p)})_{(0)}\otimes_{s_{pq}} w\cdot 1_{(2,q)}\epsilon(v\cdot 1_{(1,p)})_{(1,e)}\\
 &=&v\cdot 1_{(1,p)}\otimes_{s_{pq}}  w\cdot 1_{(2,q)}\\
 &=&v\otimes_{s_{pq}} w,
 \end{eqnarray*}
 this shows that $V\otimes_{s_{pq}} W$ satisfies the counitary condition.

 Then we check the crossed condition as follows: for all $h\in H_{pqr}, v\otimes w\in V\otimes_{s_{pq}} W$, we have
 \begin{eqnarray*}
 &&\big((v\otimes w)\cdot h_{(2,pq)}\big)_{(0)}\otimes \pi_{q^{-1}p^{-1}}(h_{(1, pqrq^{-1}p^{-1})})\big((v\otimes w)\cdot h_{(2,pq)}\big)_{(1,r)}\\
 &\stackrel{(\ref{HHQ2})}=&(v\cdot h_{(2,p)})_{(0)}\otimes (w\cdot h_{(3,q)})_{(0)}\otimes \\
 &&\pi_{q^{-1}p^{-1}}(h_{(1, pqrq^{-1}p^{-1})})\Big(\pi_{q^{-1}}((v\cdot h_{(2,p)})_{(1,qrq^{-1})}\big)(w\cdot h_{(3,q)})_{(1,r)}\Big)\\
 &\stackrel{(\ref{JR1})}=&v_{(0)}\cdot h_{(3,p)}\otimes (w\cdot h_{(5,q)})_{(0)}\otimes \pi_{q^{-1}p^{-1}}(h_{(1, pqrq^{-1}p^{-1})})\\
 &&\Big(\pi_{q^{-1}}\big(S\pi_{p^{-1}}(h_{(2,pqrq^{-1}p^{-1})})(v_{(1,qrq^{-1})}h_{(4,qrq^{-1})})\big)(w\cdot h_{(5,q)})_{(1,r)}\Big)\\
 &\stackrel{(\ref{JR1})}=&v_{(0)}\cdot h_{(3,p)}\otimes w_{(0)}\cdot h_{(6,q)}\otimes \pi_{q^{-1}p^{-1}}(h_{(1, pqrq^{-1}p^{-1})})\Big(\pi_{q^{-1}}\big(S\pi_{p^{-1}}\\ &&(h_{(2,pqr^{-1}q^{-1}p^{-1})})(v_{(1,qrq^{-1})}h_{(4,qrq^{-1})})\big)\big(S\pi_{q^{-1}}(h_{(5,qr^{-1}q^{-1})})(w_{(1,r)}h_{(7,r)})\big)\Big)\\
 &=&v_{(0)}\cdot h_{(3,p)}\otimes w_{(0)}\cdot h_{(6,q)}\otimes \pi_{q^{-1}p^{-1}}(h_{(1, pqrq^{-1}p^{-1})})\Big(\pi_{q^{-1}}\big(\pi_{p^{-1}}\\
 &&S(h_{(2,pqr^{-1}q^{-1}p^{-1})})(v_{(1,qrq^{-1})}h_{(4,qrq^{-1})})\big)\big(S\pi_{q^{-1}}(h_{(5,qr^{-1}q^{-1})})(w_{(1,r)}h_{(7,r)})\big)\Big)\\
 &=&v_{(0)}\cdot h_{(3,p)}\otimes w_{(0)}\cdot h_{(6,q)}\otimes \pi_{q^{-1}p^{-1}}(h_{(1, pqrq^{-1}p^{-1})})\Big(\big(\pi_{q^{-1}}\pi_{p^{-1}}\\
 &&S(h_{(2,pqr^{-1}q^{-1}p^{-1})})\pi_{q^{-1}}(v_{(1,qrq^{-1})}h_{(4,qrq^{-1})})\big)
 \big(S\pi_{q^{-1}}(h_{(5,qr^{-1}q^{-1})})(w_{(1,r)}h_{(7,r)})\big)\Big)\\
 &\stackrel{(\ref{SE2})}=&v_{(0)}\cdot h_{(2,p)}\otimes w_{(0)}\cdot h_{(5,q)}\otimes \pi_{q^{-1}p^{-1}}\epsilon^{t}_{pqrq^{-1}p^{-1}}(h_{(1,e)})\\
 &&\Big(\pi_{q^{-1}}(v_{(1,qrq^{-1})}h_{(3,qrq^{-1})})\big(S\pi_{q^{-1}}(h_{(4,qr^{-1}q^{-1})})(w_{(1,r)}h_{(6,r)})\big)\Big)\\
 &\stackrel{(\ref{SE5})}=&v_{(0)}\cdot 1_{(2,p)}h_{(1,p)}\otimes w_{(0)}\cdot h_{(4,q)}\otimes \big(\pi_{q^{-1}p^{-1}}S(1_{(1,pqr^{-1}q^{-1}p^{-1})})\\
 &&\pi_{q^{-1}}(v_{(1,qrq^{-1})}h_{(2,qrq^{-1})})\big)\big(S\pi_{q^{-1}}(h_{(3,qr^{-1}q^{-1})})(w_{(1,r)}h_{(5,r)})\big)\\
 &=&v_{(0)}\cdot 1_{(2,p)}h_{(1,p)}\otimes w_{(0)}\cdot h_{(4,q)}\otimes \pi_{q^{-1}}\big(S\pi_{p^{-1}}(1_{(1,pqr^{-1}q^{-1}p^{-1})})\\
 &&(v_{(1,qrq^{-1})}h_{(2,qrq^{-1})})\big)\big(S\pi_{q^{-1}}(h_{(3,qr^{-1}q^{-1})})(w_{(1,r)}h_{(5,r)})\big)\\
 &\stackrel{(\ref{JR3})}=&(v_{(0)}\cdot 1'_{(1,p)})\cdot 1_{(2,p)}h_{(1,p)}\otimes w_{(0)}\cdot h_{(4,q)}\otimes \pi_{q^{-1}}\big(S\pi_{p^{-1}}(1_{(1,pqr^{-1}q^{-1}p^{-1})})\\
 &&((v_{(1,qrq^{-1})}1'_{(2,qrq^{-1})})h_{(2,qrq^{-1})})\big)
 \big(S\pi_{q^{-1}}(h_{(3,qr^{-1}q^{-1})})(w_{(1,r)}h_{(5,r)})\big)\\
 &\stackrel{(\ref{v1h1}),(\ref{v1h3})}=&(v_{(0)}\cdot 1'_{(1,p)}1_{(2,p)})\cdot h_{(1,p)}\otimes w_{(0)}\cdot h_{(4,q)}\otimes \pi_{q^{-1}}\big(S\pi_{p^{-1}}(1_{(1,pqr^{-1}q^{-1}p^{-1})})\\
 &&((v_{(1,qrq^{-1})}1'_{(2,qrq^{-1})})h_{(2,qrq^{-1})})\big)
 \big(S\pi_{q^{-1}}(h_{(3,qr^{-1}q^{-1})})(w_{(1,r)}h_{(5,r)})\big)\\
 &=&(v_{(0)}\cdot 1_{(2,p)})\cdot h_{(1,p)}\otimes w_{(0)}\cdot h_{(4,q)}\otimes \pi_{q^{-1}}\big(S\pi_{p^{-1}}(1_{(1,pqr^{-1}q^{-1}p^{-1})})\\
 &&((v_{(1,qrq^{-1})}1_{(3,qrq^{-1})})h_{(2,qrq^{-1})})\big)
 \big(S\pi_{q^{-1}}(h_{(3,qr^{-1}q^{-1})})(w_{(1,r)}h_{(5,r)})\big)\\
 &\stackrel{(\ref{JR1})}=&(v\cdot 1)_{(0)}\cdot h_{(1,p)}\otimes w_{(0)}\cdot h_{(4,q)}\otimes \\
 &&\pi_{q^{-1}}\big((v\cdot 1)_{(1,qrq^{-1})}h_{(2,qrq^{-1})}\big)\big(S\pi_{q^{-1}}(h_{(3,qr^{-1}q^{-1})})(w_{(1,r)}h_{(5,r)})\big)\\
 &\stackrel{(\ref{JH1})}=&v_{(0)}\cdot h_{(1,p)}\otimes w_{(0)}\cdot h_{(4,q)}\otimes \\
 &&\pi_{q^{-1}}(v_{(1,qrq^{-1})})\Big(\pi_{q^{-1}}(h_{(2,qrq^{-1})})\big(S\pi_{q^{-1}}(h_{(3,qr^{-1}q^{-1})})(w_{(1,r)}h_{(5,r)})\big)\Big)\\
 &\stackrel{(\ref{SE2})}=&v_{(0)}\cdot h_{(1,p)}\otimes w_{(0)}\cdot h_{(3,q)}\otimes \\ &&\pi_{q^{-1}}(v_{(1,qrq^{-1})})\Big(\pi_{q^{-1}}\epsilon^{t}_{qrq^{-1}}(h_{(2,e)})\big(S\pi_{q^{-1}}(h_{(3,qr^{-1}q^{-1})})(w_{(1,r)}h_{(4,r)})\big)\Big)\\
 &\stackrel{(\ref{SE5})}=&v_{(0)}\cdot h_{(1,p)}\otimes w_{(0)}\cdot 1_{(2,q)}h_{(2,q)}\otimes \pi_{q^{-1}}(v_{(1,qrq^{-1})})\big(\pi_{q^{-1}}S(1_{(qr^{-1}q^{-1})})(w_{(1,r)}h_{(4,r)})\big)\\
 &\stackrel{(\ref{JR3})}=&v_{(0)}\cdot h_{(1,p)}\otimes (w_{(0)}\cdot 1'_{(1,q)})\cdot 1_{(2,q)}h_{(2,q)}\otimes \\ &&\pi_{q^{-1}}(v_{(1,qrq^{-1})})\Big(\pi_{q^{-1}}S(1_{(qr^{-1}q^{-1})})\big((w_{(1,r)}1'_{(2,r)})h_{(4,r)})\big)\Big)\\
 &=&v_{(0)}\cdot h_{(1,p)}\otimes (w_{(0)}\cdot 1_{(2,q)})\cdot h_{(2,q)}\otimes \\ &&\pi_{q^{-1}}(v_{(1,qrq^{-1})})\Big(S\pi_{q^{-1}}(1_{(qr^{-1}q^{-1})})\big((w_{(1,r)}1{(3,r)})h_{(4,r)})\big)\Big)\\
 &\stackrel{(\ref{JR1})}=&v_{(0)}\cdot h_{(1,p)}\otimes (w\cdot 1)_{(0)}\cdot h_{(2,q)}\otimes \pi_{q^{-1}}(v_{(1,qrq^{-1})})\big((w\cdot 1)_{(1,r)}h_{(3,r)}\big)\\
 &=&v_{(0)}\cdot h_{(1,p)}\otimes w_{(0)}\cdot h_{(2,q)}\otimes \pi_{q^{-1}}(v_{(1,qrq^{-1})})(w_{(1,r)}h_{(3,r)})\\
 &\stackrel{(\ref{JH2})}=&v_{(0)}\cdot h_{(1,p)}\otimes w_{(0)}\cdot h_{(2,q)}\otimes \big(\pi_{q^{-1}}(v_{(1,qrq^{-1})})w_{(1,r)}\big)h_{(3,r)}\\
 &=&(v\otimes w)_{(0)}\cdot h_{(1,pq)}\otimes (v\otimes w)_{(1,r)}h_{(2,r)},
 \end{eqnarray*}
 and
 \begin{eqnarray*}
 &&(v\otimes w)_{(0)}\cdot 1_{(1,pq)}\otimes (v\otimes w)_{(1,r)}1_{(2,r)}\\
 &\stackrel{(\ref{HHQ2})}=&(v_{(0)}\otimes w_{(0)})\cdot 1_{(1,pq)}\otimes \big(\pi_{q^{-1}}(v_{(1,qrq^{-1})})w_{((1,r))}\big)1_{(2,r)}\\
 &\stackrel{(\ref{HHQ1})}=&v_{(0)}\otimes 1_{(1,p)}\otimes w_{(0)}\cdot 1_{(2,q)}\otimes \big(\pi_{q^{-1}}(v_{(1,qrq^{-1})})w_{((1,r))}\big)1_{(3,r)}\\
 &=&v_{(0)}\otimes 1_{(1,p)}\otimes w_{(0)}\cdot 1'_{(1,q)}1_{(2,q)}\otimes \big(\pi_{q^{-1}}(v_{(1,qrq^{-1})})w_{((1,r))}\big)1'_{(2,r)}\\
 &\stackrel{(\ref{v1h1})}=&v_{(0)}\otimes 1_{(1,p)}\otimes (w_{(0)}\cdot 1'_{(1,q)})\cdot 1_{(2,q)}\otimes \big(\pi_{q^{-1}}(v_{(1,qrq^{-1})})w_{((1,r))}\big)1'_{(2,r)}\\
 &\stackrel{(\ref{JH2})}=&v_{(0)}\otimes 1_{(1,p)}\otimes (w_{(0)}\cdot 1'_{(1,q)})\cdot 1_{(2,q)}\otimes \pi_{q^{-1}}(v_{(1,qrq^{-1})})\big(w_{((1,r))}1'_{(2,r)}\big)\\
 &\stackrel{(\ref{JR3})}=&v_{(0)}\cdot 1_{(1,p)}\otimes w_{(0)}\cdot 1_{(2,q)}\otimes \pi_{q^{-1}}(v_{(1,qrq^{-1})})w_{((1,r))}\\
 &=&v_{(0)}\otimes w_{(0)}\otimes \pi_{q^{-1}}(v_{(1,qrq^{-1})})w_{((1,r))}\\
 &=&(v\otimes w)_{(0)}\otimes (v\otimes w)_{(1,r)},
 \end{eqnarray*}
 hence the crossed property of $V\otimes_{s_{pq}} W$ is hold.

 Finally, we should  check $(\ref{JH1})$ and $(\ref{JH2})$ of $V\otimes_{s_{pq}} W$ hold, and it is similar with \cite{LYZ22}.
 So we conclude that $V\otimes_{s_{pq}} W\in \mathcal Y\mathcal D\mathcal W\mathcal Q(H)_{pq}$. $\hfill \square$
 \\

 Following Turaev's left index notation, let $V\in\mathcal Y \mathcal D\mathcal W\mathcal Q(H)_{p}$, the object $^{q}V$ has the same underlying vector space as V. Given $v\in V$, we denote $^{q}v$ the corresponding element in $^{q}V$.

 \smallskip{\bf Proposition \thesection.6} Let $(V, \rho^{V})\in \mathcal Y\mathcal D\mathcal W\mathcal Q(H)_{p}$ and $q\in G$. Set $^{q}V=V$ as a vector space with structures defined by
 \begin{eqnarray}
 &&{}^{q}v\cdot h_{qpq^{-1}} ={}^{q}(v\cdot \pi_{q^{-1}}(h_{qpq^{-1}})), \label{qa}\\
 &&\rho_{r}^{^{q}V}(^{q}v)={}^{q}(v_{(0)})\otimes \pi_{q}(v_{(1,q^{-1}rq)}):= v_{<0>}\otimes v_{<1,r>}, \label{qco}
 \end{eqnarray}
 for any $v\in V$ and $h_{qpq^{-1}}\in H_{qpq^{-1}}$. Then ${}^{q}V\in \mathcal Y\mathcal D\mathcal W\mathcal Q(H)_{qpq^{-1}}$.
 \smallskip

 {\bf\emph{Proof}} We need to show that ${}^{q}V$ is a $qpq^{-1}$-Yetter-Drinfeld weak quasimodule, that is, ${}^{q}V$ is a $H_{qpq^{-1}}$-weak quasimodule, and satisfies
 coassociative, counitary and crossing conditions.
 In the following, we only show that the crossing condition holds, the others conditions are easy to check. Indeed, for all ${}^{q}v\in {}^{q}V$, we have
 \begin{eqnarray*}
 &&v_{<0>}\cdot 1_{(1,qpq^{-1})}\otimes v_{<1,r>}1_{(2,r)}\\
 &=&{}^{q}(v_{(0)})\cdot 1_{(1,qpq^{-1})}\otimes \pi_{q}(v_{1,q^{-1}rq})1_{(2,r)}\\
 &=&{}^{q}\big(v_{(0)}\cdot \pi_{q^{-1}}(1_{(1,qpq^{-1})})\big)\otimes \pi_{q}(v_{1,q^{-1}rq})1_{(2,r)}\\
 &=&{}^{q}(v_{(0)}\cdot 1_{(1,p)})\otimes \pi_{q}\big(v_{(1,q^{-1}rq)}1_{(2,r)}\big)\\
 &=&{}^{q}(v_{(0)})\otimes \pi_{q}(v_{1,q^{-1}rq})\\
 &=&v_{<0>}\otimes v_{<1,r>},
 \end{eqnarray*}
 and for all $h\in H_{qpq^{-1}r}, {}^{q}(v)\in {}^{q}(V)$, we have
 \begin{eqnarray*}
 &&({}^{q}(v)\cdot h_{(2, qpq^{-1})})_{<0>}\otimes \pi_{qp^{-1}q^{-1}}(h_{1,qpq^{-1}rqp^{-1}q^{-1}})({}^{q}(v)\cdot h_{(2, qpq^{-1})})_{<1,r>}\\
 &=&{}^{q}\big((v\cdot \pi_{q^{-1}}(h_{(2, qpq^{-1})}))_{(0)}\big)\otimes \pi_{qp^{-1}q^{-1}}(h_{1,qpq^{-1}rqp^{-1}q^{-1}})\pi_{q}\big(v\cdot \pi_{q^{-1}}(h_{(2, qpq^{-1})})_{(1,q^{-1}rq)}\big)\\
 &=&{}^{q}\big((v\cdot \pi_{q^{-1}}(h)_{(2,p)})_{(0)}\big)\otimes \pi_{q}\Big(\pi_{p^{-1}}\big(\pi_{q^{-1}}(h)_{(1,pq^{-1}rqp^{-1})}\big)\big((v\cdot \pi_{q^{-1}}(h)_{(2,p)})_{(1,q^{-1}rq)}\big)\Big)\\
 &=&{}^{q}\big((v\cdot \pi_{q^{-1}}(h)_{(1,p)})_{(0)}\big)\otimes \pi_{q}\big(v_{(1,q^{-1}rq)}\pi_{q^{-1}}(h)_{(2,q^{-1}rq)}\big)\\
 &=&{}^{q}\big((v\cdot \pi_{q^{-1}}(h_{(1,qpq^{-1})})\big)_{(0)}\big)\otimes \pi_{q}\big(v_{(1,q^{-1}rq)}\pi_{q^{-1}}(h_{(2,r)})\big)\\
 &=&{}^{q}(v_{(0)})\cdot h_{(1,qpq^{-1})}\otimes \pi_{q}(v_{(1,q^{-1}rq)})h_{(2,r)}\\
 &=&v_{<0>}\cdot h_{(1,qpq^{-1})}\otimes v_{<1,r>}h_{(2,r)}.
 \end{eqnarray*}

 This completes the proof.$\hfill \square$
 \\

 {\bf Remark} Let $(V, \rho^{V})\in \mathcal Y\mathcal D\mathcal W\mathcal Q(H)_{p}$, and let $(W, \rho^{W})\in \mathcal Y\mathcal D\mathcal W\mathcal Q(H)_{q}$, then we have ${}^{kt}V={}^{k}({}^{t}V)$ as an object in $\mathcal Y\mathcal D\mathcal W\mathcal Q(H)_{ktpt^{-1}k}$, and ${}^{k}(V\otimes_{s_{pq}} W)={}^{k}V\otimes_{s_{kpqk^{-1}}} {}^{k}W$ as an object in $\mathcal Y\mathcal D\mathcal W\mathcal Q(H)_{kpqk^{-1}}$.
 \\

 For a crossed group-cograded weak Hopf quasigroup $H$, we define $\mathcal Y\mathcal D\mathcal W\mathcal Q(H)$ as the disjoint union of all $\mathcal Y\mathcal D\mathcal W\mathcal Q(H)_{p}$ with $p\in G$.
 If we endow $\mathcal Y\mathcal D\mathcal W\mathcal Q(H)$ with tensor product as in Proposition \thesection.5, then we get the the following result.
 \\

 {\bf Theorem \thesection.7} The Yetter-Drinfeld weak quasimodule category  $\mathcal Y\mathcal D\mathcal W\mathcal Q(H)$ is a crossed category.

 {\bf\emph{Proof}} By Proposition \thesection.6 we can give a group homomorphism $\phi: G\rightarrow Aut(\mathcal Y\mathcal D\mathcal W\mathcal Q(H))$, $p\longmapsto \phi_{p}$ by
 \begin{eqnarray*}
 \phi_{p}: \mathcal Y\mathcal D\mathcal W\mathcal Q(H)_{q}\rightarrow \mathcal Y\mathcal D\mathcal W\mathcal Q(H)_{pqp^{-1}}, \qquad \phi_{p}(W)={}^{p}W,
 \end{eqnarray*}
 where the functor $\phi_{p}$ act as follows: given a morphism $f:(V,\rho^{V})\rightarrow (W,\rho^{W})$, for any $v\in V$, we set $(^{p}f)(^{p}v)={}^{p}(f(v))$.
 Then it is easy to prove $\mathcal Y\mathcal D\mathcal W\mathcal Q(H)$ is a crossed category. $\hfill \square$
 \\

 Following the ideas by $\rm{\acute{A}}$lvarez in \cite{AFG15, AFG16}, we will consider $\mathcal Y\mathcal D(H)_{p}$ the category of right-right $p$-Yetter-Drinfeld modules over $H$, which is a subcategory of $\mathcal Y\mathcal D\mathcal W\mathcal Q(H)_{p}$.
 \\

 {\bf Proposition \thesection.8} Let $(V,\rho^{V})\in \mathcal Y\mathcal D(H)_{p}$ and $(W,\rho^{W})\in \mathcal Y\mathcal D(H)_{q}$. Set ${}^{V}W={}^{p}W$ as an object in $\mathcal Y\mathcal D(H)_{pqp^{-1}}$. Define the map
 \begin{eqnarray}
 &&C_{V,W}: V\otimes W\rightarrow {}^{V}W\otimes V \nonumber \\
 &&C_{V,W}(v\otimes w)={}^{p}\big(w\cdot S^{-1}(v_{(1,q^{-1})}) \big)\otimes v_{(0)} \label{BZ}
 \end{eqnarray}
 Then $C_{V,W}$ is $H$-linear, $H$-colinear and satisfies the conditions:
 \begin{eqnarray}
 &&C_{V\otimes W, X}=(C_{V,^{W}X}\otimes id_{W})(id_{V}\otimes C_{W, X})\\
 &&C_{V, W\otimes X}=(id_{^{V}W}\otimes C_{V,X})(C_{V,W}\otimes id_{X})
 \end{eqnarray}
 for $X\in \mathcal Y\mathcal D(H)_{s}$. Moreover, $C_{^{s}V,^{s}W}={}^{s}(\cdot) C_{V,W}$.

 {\bf\emph{Proof}} Firstly, we need to show that $C_{V,W}$ is well-defined. Indeed, we have
 \begin{eqnarray*}
 &&C_{V,W}(v\cdot 1_{(1,p)}\otimes w\cdot 1_{(2,q)})\\
 &=&{}^{p}\Big((w\cdot 1_{(2,q)})\cdot S^{-1}\big((v\cdot 1_{(1,p)})_{(1,q^{-1})}\big)\Big)\otimes (v\cdot 1_{(1,p)})_{(0)}\\
 &=&{}^{p}\Big((w\cdot 1_{(4,q)})\cdot S^{-1}\big(S\pi_{p^{-1}}(1_{(1,pqp^{-1})})(v_{(1,q^{-1})}1_{(3,q^{-1})})\big)\Big)\otimes v_{(0)}\cdot 1_{(2,p)}\\
 &=&{}^{p}\Big((w\cdot 1_{(4,q)})\cdot \big(S^{-1}(1_{(3,q^{-1})})S^{-1}(v_{(1,q^{-1})})\big)S^{-1}S\pi_{p^{-1}}(1_{(1,pqp^{-1})})\big)\Big)\otimes v_{(0)}\cdot 1_{(2,p)}\\
 &=&{}^{p}\Big(w\cdot \epsilon^{t}_{q}(1_{(3,e)})S^{-1}(v_{(1,q^{-1})})\cdot S^{-1}\big(S\pi_{p^{-1}}(1_{(1,pqp^{-1})})\big)\Big)\otimes v_{(0)}\cdot 1_{(2,p)}\\
 &=&{}^{p}\Big(w\cdot S^{-1}(1'_{(2,q^{-1})})S^{-1}(v_{(1,q^{-1})})\cdot S^{-1}\big(S\pi_{p^{-1}}(1_{(1,pqp^{-1})})\big)\Big)\otimes v_{(0)}\cdot 1_{(2,p)}1'_{(1,p)}\\
 &=&{}^{p}\Big(w\cdot S^{-1}(v_{(1,q^{-1})}1_{(3,q^{-1})})\cdot S^{-1}\big(S\pi_{p^{-1}}(1_{(1,pqp^{-1})})\big)\Big)\otimes v_{(0)}\cdot 1_{(2,p)}\\
 &=&{}^{p}\Big(w\cdot S^{-1}\big(S\pi_{p^{-1}}(1_{(1,pqp^{-1})}(v_{(1,q^{-1})}1_{(3,q^{-1})}))\big)\Big)\otimes v_{(0)}\cdot 1_{(2,p)}\\
 &=&{}^{p}\big(w\cdot S^{-1})(v_{(1,q^{-1})})\big)\otimes v_{(0)}\\
 &=&C_{V,W}(v\otimes w).
 \end{eqnarray*}

 Secondly, we prove that $C_{V,W}$ is $H$-linear. For all $h\in H_{pq}$, we compute
 \begin{eqnarray*}
 &&C_{V,W}\big((v\otimes w)\cdot h_{pq}\big)\\
 &\stackrel{(\ref{HHQ1})}=&C_{V,W}(v\cdot h_{(1,p)} \otimes w\cdot h_{(2,q)})\\
 &\stackrel{(\ref{BZ})}=&{}^{p}\Big((w\cdot h_{(2,q)})\cdot S^{-1}\big((v\cdot h_{(1,p)})_{1,q^{-1}}\big)\Big)\otimes (v\cdot h_{(1,p)})_{(0)}\\
 &\stackrel{(\ref{JR2})}=&{}^{p}\Big((w\cdot h_{(4,q)})\cdot S^{-1}\big(S\pi_{p^{-1}}(h_{(1,pqp^{-1})})(v_{(1,q^{-1})}h_{(3,q^{-1})})\big)\Big)\otimes v_{(0)}\cdot h_{(2,p)}\\
 &=&{}^{p}\Big((w\cdot h_{(4,q)})\cdot S^{-1}(v_{(1,q^{-1})}h_{(3,q^{-1})})S^{-1}\big(S\pi_{p^{-1}}(h_{(1,pqp^{-1})})\big)\Big)\otimes v_{(0)}\cdot h_{(2,p)}\\
 &=&{}^{p}\Big(w\cdot h_{(4,q)})\cdot S^{-1}(h_{(3,q^{-1})})S^{-1}(v_{(1,q^{-1})})\cdot S^{-1}\big(S\pi_{p^{-1}}(h_{(1,pqp^{-1})})\big)\Big)\otimes v_{(0)}\cdot h_{(2,p)}\\
 &=&{}^{p}\Big(w\cdot \epsilon^{t}_{q}(h_{(3,e)})S^{-1}(v_{(1,q^{-1})})\cdot S^{-1}\big(S\pi_{p^{-1}}(h_{(1,pqp^{-1})})\big)\Big)\otimes v_{(0)}\cdot h_{(2,p)}\\
 &=&{}^{p}\Big(w\cdot S^{-1}(1_{(2,q^{-1})})S^{-1}(v_{(1,q^{-1})})\cdot S^{-1}\big(S\pi_{p^{-1}}(h_{(1,pqp^{-1})})\big)\Big)\otimes v_{(0)}\cdot 1_{(1,p)}h_{(2,p)}\\
 &=&{}^{p}\Big(w\cdot S^{-1}(v_{(1,q^{-1})}1_{(2,q^{-1})})S^{-1}\big(S\pi_{p^{-1}}(h_{(1,pqp^{-1})})\big)\Big)\otimes (v_{(0)}\cdot 1_{(1,p)})\cdot h_{(2,p)}\\
 &=&{}^{p}\Big(w\cdot S^{-1}(v_{(1,q^{-1})})S^{-1}\big(S\pi_{p^{-1}}(h_{(1,pqp^{-1})})\big)\Big)\otimes v_{(0)}\cdot h_{(2,p)}\\
 &=&{}^{p}\big(w\cdot S^{-1}(v_{(1,q^{-1})})\big)\cdot \pi_{p^{-1}}(h_{(1,pqp^{-1})})\otimes v_{(0)}\cdot h_{(2,p)}\\
 &=&\Big({}^{p}\big(w\cdot S^{-1}(v_{(1,q^{-1})})\big)\otimes v_{(0)}\Big)\cdot h\\
 &=&C_{V,W}(v\otimes w)\cdot h,
 \end{eqnarray*}
 so we have $C_{V,W}\big(h_{pq}\cdot (v\otimes w)\big)=h_{pq}\cdot C_{V,W}(v\otimes w)$, that is, $C_{V,W}$ is $H$-linear.

 Finally,  we prove that $C_{V,W}$ is $H$-colinear. In fact,
 \begin{eqnarray*}
 &&\rho_{r}^{^{V}W\otimes V}C_{V,W}(v\otimes w)\\
 &=&\rho_{r}^{^{V}W\otimes V}\Big(^{p}\big(w\cdot S^{-1}(v_{(1,q^{-1})}) \big)\otimes v_{(0)}\Big)\\
 &\stackrel{(\ref{HHQ2})}=&^{p}\Big(\big(w\cdot S^{-1}(v_{(1,q^{-1})}) \big)_{(0)}\Big)\otimes v_{(0)(0)}\otimes \big(w\cdot S^{-1}(v_{(1,q^{-1})})\big)_{(1,r)}v_{(0)(1,r)}\\
 &=&{}^{p}\big(w_{(0)}\cdot S^{-1}(v_{(1,q^{-1})})_{(2,q)}\big)\otimes v_{(0)(0)}\otimes \\ &&\Big(S\pi_{q^{-1}}\big(S^{-1}(v_{(1,q^{-1})})_{q,pr^{-1}p^{-1}}\big)\big(w_{(1,r)}S^{-1}(v_{(1,q^{-1})})_{(3,r)}\big)\Big)v_{(0)(1,r)}\\
 &=&{}^{p}\big(w_{(0)}\cdot S^{-1}(v_{(3,q^{-1})})\big)\otimes v_{(0)}\otimes \Big(S\pi_{q^{-1}}\big(S^{-1}(v_{(4,qrq^{-1})})\big)\big(w_{(1,r)}S^{-1}(v_{(2,r^{-1})})\big)\Big)v_{(1,r)}\\
 &=&{}^{p}\big(w_{(0)}\cdot S^{-1}(v_{(3,q^{-1})})\big)\otimes v_{(0)}\otimes S\pi_{q^{-1}}\big(S^{-1}(v_{(4,qrq^{-1})})\big)\Big(\big(w_{(1,r)}S^{-1}(v_{(2,r^{-1})})\big)v_{(1,r)}\Big)\\
 &=&{}^{p}\big(w_{(0)}\cdot S^{-1}(v_{(2,q^{-1})})\big)\otimes v_{(0)}\otimes S\pi_{q^{-1}}\big(S^{-1}(v_{(3,qrq^{-1})})\big)\big(w_{(1,r)}\tilde{\epsilon}^{s}_{r}(v_{(1,e)})\big)\\
 &=&{}^{p}\big(w_{(0)}\cdot S^{-1}(v_{(1,q^{-1})}1_{(2,q^{-1})})\big)\otimes v_{(0)}\otimes S\pi_{q^{-1}}\big(S^{-1}(v_{(2,qrq^{-1})})\big)\big(w_{(1,r)}S^{-1}(1_{(1,r^{-1}}))\big)\\
 &=&{}^{p}\big(w_{(0)}\cdot S^{-1}(v_{(1,q^{-1})}1_{(2,q^{-1})})\big)\otimes v_{(0)}\otimes \\
 &&\pi_{q^{-1}} SS^{-1}(1_{(3,qrq^{-1})})\pi_{q^{-1}}(v_{(2,qrq^{-1})})\big(w_{(1,r)}S^{-1}(1_{(1,r^{-1}}))\big)\\
 &=&{}^{p}\big(w_{(0)}\cdot S^{-1}(1)_{(2,q^{-1})}S^{-1}(v_{(1,q^{-1})}))\big)\otimes \\
 &&v_{(0)}\otimes \pi_{q^{-1}} S(S^{-1}(1)_{(1,qr^{-1}q^{-1})})\pi_{q^{-1}}(v_{(2,qrq^{-1})})\big(w_{(1,r)}S^{-1}(1)_{(3,r)}\big)\\
 &=&{}^{p}\big(w_{(0)}\cdot S^{-1}(v_{(1,q^{-1})})\big)\otimes v_{(0)}\otimes \pi_{q^{-1}}(v_{(2,qrq^{-1})})w_{(1,r)}\\
 &=&(C_{V,W}\otimes id_{H_{r}})\rho^{V\otimes W}_{r}(v\otimes w).
 \end{eqnarray*}

 This completes the proof.$\hfill \square$
 \\

 {\bf Proposition \thesection.9} Let $(V, \rho^{V})\in \mathcal Y\mathcal D(H)_{p}$ and $(W, \rho^{W})\in \mathcal Y\mathcal D(H)_{q}$. Then can give the braided $C_{V,W}$ an inverse $C^{-1}_{V, W}$, which is defined by
 \begin{eqnarray*}
 C^{-1}_{V, W}: {}^{V}W\otimes V&\rightarrow& V\otimes W,\\
 C^{-1}_{V, W}({^{p}w\otimes v})&=&v_{(0)}\otimes w\cdot v_{(1,q)} ,
 \end{eqnarray*}
 where $p,q\in G$.

 {\bf\emph{Proof}} For any $v\in V, w\in W$, we have
 \begin{eqnarray*}
 C^{-1}_{V, W}C_{V,W}(v\otimes w)
 &=&C^{-1}_{V, W}({}^{p}(w\cdot S^{-1}(v_{(1,q^{-1})}))\otimes v_{(0)})\\
 &=&v_{(0)(0)}\otimes (w\cdot S^{-1}(v_{(2,q^{-1})}))\cdot v_{(0)(1,q)}\\
 &=&v_{(0)}\otimes (w\cdot S^{-1}(v_{(2,q^{-1})}))\cdot v_{(1,q)}\\
 &=&v_{(0)}\otimes w\cdot S^{-1}\epsilon^{s}_{q^{-1}}(v_{(1,e)})\\
 &=&v_{(0)}\cdot 1_{(1,p)}\otimes w\cdot \tilde{\epsilon}^{s}_{q}(v_{(1,e)}1_{(2,e)})\\
 &=&v_{(0)}\cdot 1_{(1,p)}\otimes w\cdot \tilde{\epsilon}^{s}_{q}\big(v_{(1,e)}\tilde{\epsilon}^{s}_{e}(1_{(2,e)})\big)\\
 &=&v_{(0)}\cdot 1_{(1,p)}\otimes w\cdot \tilde{\epsilon}^{s}_{q}(v_{(1,e)})\tilde{\epsilon}^{s}_{q}(1_{(2,e)})\\
 &=&v_{(0)}\cdot 1_{(1,p)}\otimes w\cdot \tilde{\epsilon}^{s}_{q}(v_{(1,e)})1_{(2,q)}\\
 &=&v_{(0)}\cdot \epsilon^{s}_{p}(v_{(1,e)})1_{(1,p)}\otimes w\cdot1_{(2,q)}\\
 &=&v_{(0)}\cdot 1'_{(1,p)}\epsilon(v_{(1,e)}1'_{(2,e)})1_{(1,p)}\otimes w\cdot1_{(2,q)}\\
 &=&v_{(0)}\epsilon(v_{(1,e)})\cdot 1_{(1,p)}\otimes w\cdot1_{(2,q)}\\
 &=&v\cdot 1_{(1,p)}\otimes w\cdot1_{(2,q)}\\
 &=&v\otimes w.
 \end{eqnarray*}
 Conversely, for any ${}^{p}w\in {}^{V}W, v\in V$,
 \begin{eqnarray*}
 C_{V, W}C^{-1}_{V,W}({}^{p}w\otimes v)
 &=&C_{V, W}(v_{(0)}\otimes w\cdot v_{(1,q)})\\
 &=&{}^{p}\Big(w\cdot \big(v_{(2,q)}\cdot S^{-1}(v_{(1,q^{-1})})\big)\Big)\otimes v_{(0)}\\
 &=&{}^{p}\big(w\cdot v_{(2,q)}S^{-1}(v_{(1,q^{-1})})\big)\cdot 1_{(1,pqp^{-1})}\otimes v_{(0)}\cdot 1_{(2,p)}\\
 &=&{}^{p}\big(w\cdot v_{(2,q)}S^{-1}(v_{(1,q^{-1})})\pi_{p^{-1}}(1_{(1,pqp^{-1})})\big)\otimes v_{(0)}\cdot 1_{(2,p)}\\
 &=&{}^{p}\Big(w\cdot S^{-1}\big(v_{(1,q^{-1})}S(v_{(2,q)})\big)\pi_{p^{-1}}(1_{(1,pqp^{-1})})\Big)\otimes v_{(0)}\cdot 1_{(2,p)}\\
 &=&{}^{p}\big(w\cdot S^{-1}\epsilon^{t}_{q^{-1}}(v_{(1,e)})\pi_{p^{-1}}(1_{(1,pqp^{-1})})\big)\otimes v_{(0)}\cdot 1_{(2,p)}\\
 &=&{}^{p}\big(w\cdot \tilde{\epsilon}^{t}_{q^{-1}}(v_{(1,e)})\pi_{p^{-1}}(1_{(1,pqp^{-1})})\big)\otimes v_{(0)}\cdot 1_{(2,p)}\\
 &=&{}^{p}\big(w\cdot \pi_{p^{-1}}(\tilde{\epsilon}^{t}_{pqp^{-1}}\pi_{p}(v_{(1,e)})1_{(1,pqp^{-1})})\big)\otimes v_{(0)}\cdot 1_{(2,p)}\\
 &=&{}^{p}\big(w\cdot \pi_{p^{-1}}(1_{(1,pqp^{-1})})\big)\otimes v_{(0)}\cdot \epsilon^{s}_{p}\pi_{p}(v_{(1,e)})1_{(2,p)}\\
 &=&{}^{p}\big(w\cdot \pi_{p^{-1}}(1_{(1,pqp^{-1})})\big)\otimes v_{(0)}\cdot 1'_{(1,p)}\epsilon\big(\pi_{p}(v_{(1,e)})1'_{(2,e)}\big)1_{(2,p)}\\
 &=&{}^{p}\big(w\cdot \pi_{p^{-1}}(1_{(1,pqp^{-1})})\big)\otimes v_{(0)}\epsilon\big(\pi_{p}(v_{(1,e)})\big)\cdot 1_{(2,p)}\\
 &=&{}^{p}\big(w\cdot \pi_{p^{-1}}(1_{(1,pqp^{-1})})\big)\otimes v_{(0)}\epsilon(v_{(1,e)})\cdot 1_{(2,p)}\\
 &=&{}^{p}\big(w\cdot \pi_{p^{-1}}(1_{(1,pqp^{-1})})\big)\otimes v\cdot 1_{(2,p)}\\
 &=&{}^{p}w\cdot 1_{(1,pqp^{-1})}\otimes v\cdot 1_{(2,p)}\\
 &=&{}^{p}w\otimes v.
 \end{eqnarray*}
 Since then $C_{V,W}$ is an isomorphism with inverse $C^{-1}_{V, W}$.$\hfill \square$
 \\


 As a consequence of the above results, we obtain another main result of this paper.
 \\

 {\bf Theorem \thesection.10} For a crossed group-cograded weak Hopf quasigroup $H$, we define $\mathcal Y\mathcal D(H)$ as the disjoint union of all $\mathcal Y\mathcal D(H)_{p}$ with $p\in G$. Then $\mathcal Y\mathcal D(H)$ is a braided crossed category over group $G$.

 {\bf\emph{Proof}} As for $\mathcal Y\mathcal D(H)$ is a subcategory of the category $\mathcal Y\mathcal D\mathcal W\mathcal Q(H)$, so it is a crossed category.
 Then we only need prove $\mathcal Y\mathcal D(H)$ is braided.

 The braiding in $\mathcal Y\mathcal D(H)$ can be given by Proposition \thesection.8, and the braiding is invertible, its inverse is the family $C^{-1}_{V,W}$, which defined in Proposition \thesection.9.
 Hence it is obvious that $\mathcal Y\mathcal D(H)$ is a braided crossed category. $\hfill \square$

\section*{Acknowledgements}

 The work was partially funded by China Postdoctoral Science Foundation (Grant No.2019M651764), National Natural Science Foundation of China (Grant No. 11601231).

\end {document}